\documentclass{article}
\usepackage{slashbox}
\usepackage[T1]{fontenc}
\usepackage{amsmath}
\usepackage{amssymb}
\usepackage{color}
\usepackage{graphicx}
\usepackage{rotating}
\usepackage{algorithmic} 
\usepackage[lined,boxed,commentsnumbered]{algorithm2e}
\usepackage{lineno,hyperref}
\newcommand{\bi}{ \boldsymbol{i}}

\newcommand{\K}{ \mathcal{K} }
\newcommand{\A}{ \mathcal{A} }
\newcommand{\M}{ \mathcal{M} }
\newcommand{\Q}{ \mathbb{Q} }
\newcommand{\W}{ \mathrm{W} }
\newcommand{\BB}{\mathrm{B}}
\newcommand{\F}{ \boldsymbol{F} }

\newcommand{\I}{ \boldsymbol{i} }
\newcommand{\XI}{ \boldsymbol{\xi} }
\newcommand{\q}{ \boldsymbol{q} }
\newcommand{\x}{ \boldsymbol{x} }
\newcommand{\uu}{ \boldsymbol{u} }
\newcommand{\f}{ \boldsymbol{f} }

\newtheorem{remark}{Remark}

\def\GS{\color{black}} 
\def\B{\color{black}}
\def\MT{\color{black}}

\begin{document}
\bibliographystyle{siam}

\pagestyle{myheadings}
\markboth{G. Sangalli and M. Tani}{Matrix-free isogeometric analysis}

\title { \GS Matrix-free  weighted quadrature for a computationally efficient isogeometric  $k$-method}

\author{Giancarlo Sangalli$^\dagger$
\and  Mattia Tani\thanks{%
Universit\`a di Pavia, Dipartimento di Matematica ``F. Casorati'', 
Via A. Ferrata 1, 27100 Pavia, Italy. Emails: 
{\tt  \{sangia05,mattia.tani\}@unipv.it}}}
%
\maketitle
\begin{abstract} 

\GS 
The $k$-method is the isogeometric method based on  splines (or NURBS, etc.) with maximum regularity.  
{\MT When implemented following the paradigms of classical finite element methods, the computational resources required by the $k-$method are prohibitive even for moderate degree.}
In order to address this issue, we propose  a matrix-free strategy  combined with   weighted quadrature, which is an ad-hoc strategy to compute the integrals of the Galerkin system. Matrix-free  weighted quadrature (MF-WQ) speeds up matrix operations, and, perhaps even more important, greatly reduces memory consumption.  Our strategy  also requires  an efficient preconditioner for the linear system iterative solver. In this work we deal with an elliptic model problem, and adopt a  preconditioner based on the Fast Diagonalization method, an old idea to solve Sylvester-like equations. Our numerical tests show that the isogeometric solver based on  MF-WQ is faster than standard approaches (where the main cost is the matrix  formation  by standard Gaussian quadrature) even for low degree. But the main achievement  is that, with  MF-WQ, the  $k$-method gets orders of magnitude faster by increasing the degree, given a target accuracy.
Therefore, we are  able to show the superiority,  in terms of computational efficiency, of the high-degree $k$-method with respect to  low-degree isogeometric  discretizations.  What we present here  is applicable to more complex and realistic differential problems,  but its effectiveness will depend on the preconditioner stage, which is  as always problem-dependent. This situation is typical of modern high-order methods: the overall performance is mainly related to the quality of the preconditioner.

\vskip 1mm
\noindent
{\bf Keywords:} Isogeometric analysis, $k$-method,  matrix-free, weighted quadrature, preconditioner.
\end{abstract}

\section{Introduction}

Introduced in the seminal paper \cite{Hughes2005}, isogeometric analysis  is a numerical method to solve partial differential equations (PDEs).
 It can be seen as an extension of the standard finite element method, where  the unknown solution of the PDE is approximated by the same functions that are adopted in computer-aided design for the parametrization of  the PDE domain,  typically splines and extensions, such as non-uniform rational B-splines (NURBS). We refer to the monograph \cite{Cottrell2009} for a detailed description of this approach.

One feature that distinguishes isogeometric analysis  from finite element analysis  is the regularity of the basis functions. Indeed, while finite element functions are typically $C^0$, the global regularity of  splines of degree $p$ goes up to $C^{p-1}$. We denote   $k$-method the isogeometric method where splines of highest possible regularity are employed. In this setting, higher accuracy is achieved by  $k$-refinement, that is, raising the degree $p$ and refining the mesh, with $C^{p-1}$ regularity at the new inserted knots, see \cite{Hughes2005}.
  
The  $k$-method  leads to several advantages, such as  higher accuracy per degree-of-freedom and improved spectral accuracy,  see for example  \cite{Hughes2005,Cottrell2007,BeiraodaVeiga2011,BeiraodaVeiga2014,Hughes2008,evans2009n}. {\MT At the same time, the  $k$-method brings significant challenges at the computational level. Indeed, using standard  finite element routines, its computational cost  grows too fast with respect to the degree $p$, making degree raising unfavourable from the viewpoint of computational efficiency. As a consequence, in complex isogeometric  simulations quadratic $C^1$ NURBS are preferred, see \cite{Bazilevs2007,morganti2015patient}.}

The computational cost of the $k$-method  is well understood in literature. The standard formation of the Poisson (Galerkin) stiffness matrix, by Gauss quadrature and element assembly, takes $O(N p^{9})$ floating point operations (flops) in 3D,  where  $N$ is the number of degrees of freedom (see, e.g., \cite{Calabro2017} for the details). Computational limitations appear also in the solution phase of the linear system: the large support overlap associated with B-splines makes the associated Galerkin matrices less sparse than their finite element counterparts, which decreases the performance of direct solvers:  it has been shown in \cite{Collier2012} that a multifrontal direct method would require $O(N^2 p^3)$ flops and $O(N^{4/3} p^2)$ memory to solve the system.

Several alternatives to reduce this computational effort have been proposed in recent years. For what concerns the formation of isogeometric matrices, we refer to  \cite{Antolin2015,Calabro2017,Mantzaflaris2017,Hughes2010,Johannessen2016,Barton2016a,Schillinger2014,hiemstra2017optimal}. On the linear algebra side, the focus of research shifted to the development of preconditioned iterative approaches, see \cite{Buffa2013,Gahalaut2013,Donatelli2015,Hofreither2016,Sangalli2016,Tani2017}.

The starting points of this paper are the results of \GS \cite{Calabro2017}, where a new strategy, named weighted quadrature, is proposed  to form isogeometric Galerkin matrices. \B
With weighted quadrature, integrals are written  by incorporating the test function in the integral measure, then  defining different  quadrature weights for each row of the matrix. The advantage of weighted quadrature is that, for the $k$-method, the number of quadrature points is  essentially independent of $p$, and significantly less  than for other quadrature approaches. Exploiting the tensor-structure of multivariate splines, weighted quadrature allows to  form the stiffness matrix for  the 3D Poisson problem in  $O(N p^{4})$ flops,   close to the ideal formation cost of $O(N p^3)$ flops, which is the number of nonzero entries of the matrix. 

\GS  In order to further improve the computational efficiency of the $k$-method,  we now  abandon the idea of forming and storing the stiffness and mass matrices and, still relying  on weighted quadrature as in  \cite{Calabro2017},  consider a matrix-free approach.  Therefore  In such a case, the system matrix is available only as a function that computes matrix-vector products.  This is exactly what is needed by an iterative solver.
Matrix-free approaches have been use in  high-order methods based on a tensor construction like spectral elements (see \cite{tufo1999terascale})  and have been recently extended to  $hp$-finite elements \cite{Kirby2011,AINSWORTH2016140}. They are commonly used in non-linear solvers, parallel implementations, typically  for application that are computationally  demanding, for example in computational-fluid-dynamics \cite{JOHAN1991281,RASETARINERA2001718}.\B

The cost to initialize our matrix-free weighted quadrature (MF-WQ) approach is only $O(N)$ flops, while the computation of matrix-vector products costs only $O(N p)$ flops. Moreover,  the memory required by MF-WQ is just $O(N)$, i.e., it is proportional to the number of degrees of freedom. On the other hand, in 3D problems  the memory required to store the matrix would be $O(N p^3)$, 
and the cost to compute standard matrix-vector products would be $O(N p^3)$ flops. It is important to remark that, while in some cases the reduction in storage is the major motivation of matrix-free approaches,  in our framework both flops and memory savings are fundamental in order to make the use of the  high-degree $k$-method possible and advantageous. Indeed, as will be pointed out, other matrix-free approaches which rely on more standard quadrature rules (e.g. Gaussian quadrature) require $O(N p^{4})$ flops to compute a matrix-vector product. {\MT Note that this is even more costly than a standard product}.

\GS Clearly, with iterative solvers it is essential to have a good preconditioner at disposal. For this reason, we consider in here  an elliptic model problem, for which  we have developed in  \cite{Sangalli2016}  a preconditioner that, again, takes advantage of the  tensor-product structure of multivariate splines. Each application of the preconditioner is a fast direct solve of a   suitable  Sylvester-like equation. The  preconditioner is   robust with respect to $p$ and to the mesh size. Its computational cost is independent of $p$, and numerical evidence from  \cite{Sangalli2016}  shows that such cost is negligible with respect to the matrix formation, storage,  and matrix-vector operations. In \cite{Tani2017}, this preconditioner was discussed in the context of weighted quadrature, where the matrices involved are slightly non-symmetric, showing that the same properties apply in this setting. 
\B

We combine MF-WQ  and  the preconditioner of \cite{Sangalli2016} in an innovative approach which is, in the case of the $k$-method, orders of magnitude faster than the standard approach inherited by finite elements. The speedup on a mesh of  $256^3$ elements is 13 times for degree $p=1$, 44 times for degree $p=2$, while higher degrees can not be handled in the standard framework.  Indeed, in the standard approach, higher degrees are beyond the memory constraints of our workstation, while they are easily allowed in the new framework. In one of our tests,  we approximate  a solution with   oscillating behaviour  with   target accuracy of  around $10^{-3}$ (relative) error in $H^1$ norm,  and get the following computing times

  \begin{itemize}
  \item {\MT $3.0 \cdot 10^4$}  seconds  for  quadratic $C^1$ approximation on a  mesh of $256^3$ elements, with the  standard approach (system matrix formed by Gaussian integration, conjugate gradient solver with preconditioner from \cite{Sangalli2016}),
 \item {\MT $6.9 \cdot 10^2$} seconds  for quadratic $C^1$ approximation on a  mesh of $256^3$ elements, with the new approach (MF-WQ, BiCGStab iterative solver with preconditioner from \cite{Sangalli2016}),
  \item {\MT $2.7 \cdot 10^0$} seconds for  degree $p=8$,  $C^7$ spline  approximation  on a  mesh of $32^3$ elements, with the new approach (as above).
  \end{itemize}
\GS The numerical tests are performed in MATLAB. This likely favors the standard approach since its dominant cost is related  to  the matrix-vector multiplication, which is performed by a (fast) compiled routine.  However, considering only the new approach, we notice a speedup factor higher than $250$ by raising the degree from $p=2$ to $p=8$:  eventually, this gives clear evidence of the superiority of the high-degree $k$-method with respect to low-degree isogeometric discretizations in terms of computational efficiency.
\B

\GS MF-WQ has been also  implemented and tested in an innovative environment and hardware for dataflow computing, in the thesis \cite{tesi-Ruben}. 
For more complex and realistic differential problems, the overall approach \B is still applicable but its effectiveness depends on the preconditioning strategy, which is as always problem dependent. This situation is typical of modern high-order methods: the final performance is mainly related to the quality of the preconditioner. We have already obtained promising results for the Stokes system \cite{montardini2017robust}, and will devote further studies to the topic.

The outline of the paper is as follows. In Section \ref{sec:preliminaries} we introduce the Poisson model problem, its isogeometric discretization and briefly describe the weighted quadrature approach. In Section \ref{sec:matrix-free} we extend it  to the  matrix-free framework. The whole algorithm design is summarized  in Section \ref{sec:experiments} where also numerical experiments are reported. Finally, in Section \ref{sec:conclusions} we draw the conclusions.

\section{Preliminaries} 
\label{sec:preliminaries}

In this section we present the fundamental ideas and notation on isogeometric analysis, weighted quadrature and  Kronecker matrices.
\subsection{Splines-based isogeometric method}
\label{sec:splin-isog-meth}

{\MT The aim of this work is to discuss MF-WQ, focusing on the action of the  mass and the stiffness operators arising in Galerkin discretizations of elliptic equations. For this reason we consider the model problem: 
\begin{equation}\label{eq:poisson}
\left\{
\begin{aligned}
- \text{div} \left( {K} \nabla \mathrm{u} \right) + \alpha \mathrm{u} = f &\quad \text{on} \quad \Omega, \\
\mathrm{u} = 0 & \quad \text{on} \quad \partial \Omega,
\end{aligned}
\right.
\end{equation}
where $\Omega$ is a domain of $\mathbb{R}^d$ (mainly, we will be interested in the case $d=3$), $K$ is a symmetric and positive definite $d \times d$ matrix and $\alpha \geq 0$. We assume that both $K$ and $\alpha$ are smooth functions.
}

Given two positive integers $p$ and $m$, for each direction $l=1,\ldots,d$ we consider the open knot vector $\Xi_l = \left\{\xi_1^{(l)}, \ldots, \xi^{(l)}_{m+p+1} \right\}$, with
$$ 0 = \xi^{(l)}_1 =\ldots=\xi^{(l)}_{p+1} < \xi^{(l)}_{p+2} \le \ldots \le \xi^{(l)}_{m} < \xi^{(l)}_{m+1}=\ldots=\xi^{(l)}_{m+p+1} = 1, $$
where, besides the first and the last, knots can be repeated up to multiplicity $p$.


We denote with $\hat{b}_{l,i}$, $i = 0, 1, \ldots, m - 1$, the univariate B-splines (the spline  basis functions) of degree $p$ corresponding to the knot vector $\Xi_l$, $l = 1, \ldots,d$ (we refer to \cite[Chapter 2]{Cottrell2009} for the details). These are piecewise polynomials of degree $p$ which are $C^{p-r}$ continuous at knots with multiplicity $r$.
Since we are mainly interested in the $k-$method, we restrict to the case of $C^{p-1}$ regularity, i.e. all knots have multiplicity 1 except the first and the last. 
The intervals $\left[\xi^{(l)}_{k},\xi^{(l)}_{k+1} \right]$, $k = p+1,\ldots,m$, are are referred to as \textit{knot-spans}.

Multivariate B-splines are defined as tensor products of univariate B-splines. Given a multi-index $\I = (i_1, \ldots, i_d)$, we define
\begin{equation} \label{eq:multivariatesplines} \hat{b}_{\I}(\XI) = \hat{b}_{1,i_1}\left(\xi_1\right) \ldots \hat{b}_{d,i_d} \left(\xi_d\right), \end{equation}
where $\XI = \left(\xi_1,\ldots,\xi_d\right)$. To avoid proliferation of notation, we assume that all the univariate spline spaces have the same degree and the same dimension. We emphasize, however, that this does not restrict the approach described in this work.

In  isogeometric analysis, $\Omega$ is typically  given by a NURBS
parametrization. Just for the sake of simplicity, we consider a
single-patch spline parametrization though everything we say works for NURBS as well (indeed we will use NURBS in the  numerical benchmarks of Section \ref{sec:experiments}).
\begin{displaymath}
  \Omega = \boldsymbol{F} ([0,1]^d), \text{ with }
  \boldsymbol{F}(\boldsymbol \xi) = \sum_{\bi}
  \boldsymbol{C}_{\bi} \hat  b_{\bi} (\boldsymbol \xi),
\end{displaymath}
where $  \boldsymbol{C}_{\bi} $ are the control points. Following the 
isoparametric paradigm, the basis functions $b_{\bi}$ on  $\Omega$
are defined as $b_{\bi} = \hat b_{\bi}\circ \boldsymbol{F} ^{-1}$.  
We define the index set
$$ \mathcal{I} = \left\{ (i_1,\ldots,i_d) \in \mathbb{N}^d \left.\right| \; 1 \leq i_l \leq m-2, \; 1\leq l \leq d \right\}. $$
The isogeometric space, incorporating the homogeneous Dirichlet
boundary condition, then reads
\begin{equation} \label{eq:space} V_h = \left\{ \hat{b}_{\I} \circ \F^{-1} \left.\right| \I \in \mathcal{I} \right\}, \end{equation}
If we define $n := m - 2$, the total number of degrees of freedom is then $N := \mbox{dim} (V_h) = n^d$. 
In the following, we will identify the multi-index $\I = (i_1,\ldots,i_d)$ with the scalar index $i = 1+ \sum_{l=1}^d n^{l-1} (i_l - 1)$.

Given a matrix $A$, we denote with $A_{ij}$ its $(i,j)$ entry.
Then the Galerkin mass matrix $\M$ is defined as
\begin{equation}
\label{eq:M}
  \begin{aligned}
    \M_{ij} & =   \int_{\Omega} {\MT \alpha(\boldsymbol{x})} b_{i}(\boldsymbol{x}) b_{j}(\boldsymbol{x}) d\boldsymbol{x} \\
& =   \int_{[0,1]^d} c\left(\boldsymbol{\xi}\right) \hat{b}_{i}\left(\boldsymbol{\xi}\right) \hat{b}_j \left(\boldsymbol{\xi}\right)\; d\boldsymbol{\xi}, \qquad i,j = 1,\ldots,N,
  \end{aligned}
\end{equation}
where $c\left(\boldsymbol{\xi}\right) = {\MT \alpha\left(\boldsymbol{\xi}\right) }\mbox{det}\left(J_{\F}(\boldsymbol{\xi})\right)$ and $J_{\F}$ denotes the Jacobian of $\F$. Similarly, the Galerkin stiffness matrix $\mathcal{K}$ is defined as
\begin{equation}
\label{eq:A}
  \begin{aligned}
    \K_{ij} & =   \int_{\Omega} \left(\nabla b_{i}(\boldsymbol{x}) \right)^T {\MT K (\boldsymbol{x})}
  \nabla b_{{j}}(\boldsymbol{x}) d\boldsymbol{x}\\
& =   \int_{[0,1]^d} \left(\nabla \hat{b}_i(\boldsymbol{\xi})\right)^T
C(\boldsymbol{\xi}) \,  \nabla \hat{b}_j \left(\boldsymbol{\xi}\right)\; d\boldsymbol{\xi}
, \qquad i,j = 1,\ldots,N,
  \end{aligned}
\end{equation}
where
\begin{equation} \label{eq:C} C = \mbox{det}\left(J_{\F}\right)  J_{\F}^{-T} {\MT K (\boldsymbol{x})} 
J_{\F}^{-1}. \end{equation} 

The discrete version of problem \eqref{eq:poisson} then reads
\begin{equation} \label{eq:linearsystem} \left(\K + \M \right)\uu =  \A \uu = \f, \end{equation}
where
$$\f_i = \int_{[0,1]^d} c\left(\boldsymbol{\xi}\right) \hat{b}_{i}(\boldsymbol{\xi}) f(\boldsymbol{\xi}) d \boldsymbol{\xi} $$
and $\uu \in \mathbb{R}^N$ is the vector of coefficients of the Galerkin solution.

\subsection{Weighted quadrature}
\label{sec:weighted-quadrature}

Weighted quadrature was first introduced in \cite{Calabro2017} to efficiently compute Galerkin matrices in isogeometric analysis. Here we describe its main features, and in particular those that are relevant for MF-WQ, and refer to \cite{Calabro2017} for further details.

We start by discussing the case of the mass matrix. The $(i,j)$ entry of $\M$ is approximated using a special quadrature rule $\Q_{i}$ that incorporates the trial function $\hat{b}_{i}$ and treats $c \hat{b}_{j}$ as the integrand function. More precisely,
\begin{equation} \label{eq:Mwq} \M_{ij} \approx \widetilde{\M}_{ij} = \Q_{i}\left( c\left(\cdot\right) \hat{b}_{j} \left(\cdot\right)\right) = \sum_{\q \in \mathcal{Q}} w_{i,\q} c(\x_{\q}) \hat{b}_{j}(\x_{\q}).\end{equation}
Here $\mathcal{Q}$ is the multi-index set
$$ \mathcal{Q} =  \left\{ (q_1,\ldots,q_d) \in \mathbb{N}^n\left.\right| \; 1 \leq q_l \leq n_q, \; 1\leq l \leq d \right\}, $$
the tensor set of quadrature points is $\left\{\x_{\q}= \left(x_{1,q_1},\ldots,x_{d,q_d}\right) \right\}_{\q \in \mathcal{Q}}$ and the quadrature weights, related to $\Q_{i} $, are $\left\{ w_{i,q} \right\}_{\q \in \mathcal{Q}}$. 
The total number of quadrature points/weights in \eqref{eq:Mwq}  is denoted with $N_q := n_q^d$. Similarly as before, we identify the multi-index $\q = (q_1,\ldots,q_d)$ with the scalar index $q = 1+ \sum_{l=1}^d n_q^{l-1} (q_l - 1) $.

The quadrature weights are chosen such that the following exactness conditions are satisfied:
\begin{equation} \label{eq:exactness_M}
\int_{[0,1]^d} \hat{b}_{i}\left(\boldsymbol{\xi}\right) \hat{b}_j \left(\boldsymbol{\xi}\right)\; d\boldsymbol{\xi} = \Q_{i}\left( \hat{b}_{j} \left(\cdot\right)\right) = \sum_{q = 1}^{N_q} w_{i,q} \hat{b}_{j}(\x_{q}), \qquad j = 1,\ldots,N;
\end{equation}  
instead, the quadrature points  are selected a priori, such that \eqref{eq:exactness_M} is a well posed problem for the quadrature weights (see below).  Moreover the  quadrature points  are independent on the index $i$. This is the framework developed in \cite{Calabro2017}, but there are other possibilities, such as the  Gaussian weighted quadrature, currently under study.

Because of the tensor structure of the B-splines and of the quadrature points, the quadrature  weights can be determined by tensorization.
For each direction $l = 1,\ldots,d$, and for each $i_l = 1,\ldots,n$, the set $\left\{w_{l,i_l,q_l}\right\}_{q_l = 1,\ldots,n_q}$ is chosen to satisfy the univariate analogous of the exactness conditions \eqref{eq:exactness_M}, i.e.
\begin{equation} \label{eq:exactness1d} 
\int_0^1 \hat{b}_{l,i_l}\left(\xi\right) \hat{b}_{l,j_l} \left(\xi\right)\; d \xi = \sum_{q_l = 1}^{n_q} w_{l,i_l,q_l} \hat{b}_{l,j_l}(x_{l,q_l}), \qquad j_l = 1,\ldots,n. 
\end{equation}
Moreover, since the univariate quadrature weights $w_{i_l,q_l}$, $q_l=1,\ldots,n_q$, approximate the integral measure associated with the basis function $\hat{b}_{l,i_l}$, we also impose that
\begin{equation} \label{eq:weights-support} w_{i_l,q_l} = 0 \qquad \mbox{if} \qquad x_{l,q_l} \notin \mbox{supp} \left(\hat{b}_{l,i_l}\right) .
\end{equation}
It is shown in \cite{Calabro2017} that, in the case of maximum regularity, the conditions \eqref{eq:exactness1d} and  \eqref{eq:weights-support} are well posed by taking 
as quadrature points the midpoints and the endpoints of each knot span 
(expect for the first and the last, where a few more points are needed). 

There are  two consequences, that are important for the computational aspects of MF-WQ. First, since for a given $l=1,\ldots,d$ and $i_l = 1,\ldots,n$ the support of each univariate basis functions $\hat{b}_{l,i_l}$ covers at most $p+1$ knot spans, there are only roughly $2p = O(p)$ quadrature weights $w_{l,i_l,q_l}$, $q_l = 1,\ldots,n_q$ that are nonzero.
Second, the number of quadrature points is about twice the number of knot spans, i.e. $n_q \approx 2 n$. In particular, the number of quadrature points does not depend on $p$. This represents a major improvement with respect to traditional Gaussian-like quadrature approaches where the total  number of quadrature points typically grows as $O(n p)$.

The quadrature weights that satisfy the multivariate exactness conditions \eqref{eq:exactness_M} are then given by:
\begin{equation} 
w_{i,q} = w_{1,i_1,q_1} \ldots w_{d,i_d,q_d}.
\end{equation}
%

As discussed in \cite{Calabro2017}, when weighted quadrature is coupled with sum-factorization (which exploits the tensor structure to save computations), the total cost for forming the mass matrix is $O(N p^{d+1})$ flops. 

The stiffness matrix can be calculated similarly. 
We first write the $(i,j)$ entry of $\K$ from \eqref{eq:A} as 
\begin{equation} \label{eq:stiffness_sum}
 \K_{ij} =  \sum_{\alpha,\beta = 1}^d \int_{[0,1]^d} C_{\alpha \beta}(\boldsymbol{\xi}) \, \partial_{\alpha} \hat{b}_i(\boldsymbol{\xi}) \, \partial_{\beta} \hat{b}_j \left(\boldsymbol{\xi}\right)\; d\boldsymbol{\xi}  
\end{equation}
Each integral in this sum is approximated separately using a specific quadrature rule, incorporating $\partial_{\alpha} \hat{b}_i$ in the integral measure and treating $C_{\alpha \beta} \, \partial_{\beta} \hat{b}_j$ as the integrand function. Thus, we have
\begin{equation} \label{eq:Kwq}
\K_{ij} \approx \widetilde{\K}_{ij} = \Q_{i}^{(\alpha, \beta)}\left( C_{\alpha \beta}\left(\cdot\right) \partial_{\beta} \hat{b}_{j} \left(\cdot\right)\right) = \sum_{q = 1}^{N_q} w_{i,q}^{(\alpha,\beta)} C_{\alpha \beta}(\x_{q}) \; \partial_{\beta} \hat{b}_{j}(\x_{q}).
\end{equation}
We emphasize that here the quadrature weights depend not only on $i$ but also on $\alpha$ and $\beta$. Following \cite{Calabro2017}, they are chosen to satisfy the exactness conditions
\begin{equation} \label{eq:exactness_K}
\int_{[0,1]^d} \partial_{\alpha} \hat{b}_{i}\left(\boldsymbol{\xi}\right) \partial_{\beta} \hat{b}_j \left(\boldsymbol{\xi}\right)\; d\boldsymbol{\xi} = \Q_{i}^{(\alpha,\beta)}\left( \partial_{\beta} \hat{b}_{j} \left(\cdot\right)\right) := \sum_{q = 1}^{N_q} w^{(\alpha,\beta)}_{i,q} \partial_{\beta} \hat{b}_{j}(\x_{q}),
\end{equation}
for every $j = 1,\ldots,N$.  Exactness conditions different than \eqref{eq:exactness_K} can be considered. Finally, the weights (and points) can be constructed by tensorization
\begin{equation}
w_{i,q}^{(\alpha,\beta)} = w^{(\alpha,\beta)}_{1,i_1,q_1} \ldots w^{(\alpha,\beta)}_{d,i_d,q_d}.
\end{equation}
We remark that the multi-index set $\mathcal{Q}$, and the set of points $\left\{x_q\right\}_{q \in \mathcal{Q}}$ can be chosen exactly as in the mass case. Thus, the cost for forming the stiffness matrix is again $O(N p^{d+1})$ flops, although it is $d^2$ more expensive than for the mass matrix, due to the $d^2$ integrals appearing in \eqref{eq:stiffness_sum}.

\subsection{Kronecker product and sum-factorization}

Let $A, B \in \mathbb{R}^{ s \times t}$. The \textit{Kronecker product} between $A$ and $B$ is defined as
$$ A \otimes B = \begin{bmatrix} A_{11} B & \ldots & A_{1 t} B \\ \vdots & \ddots & \vdots \\ A_{s 1} B & \ldots & A_{s t} B \end{bmatrix} \; \in \mathbb{R}^{ s^2 \times t^2}, $$
The Kronecker product is an associative operation, and it is bilinear with respect to matrix sum and scalar multiplication. 
We refer to \cite{VanLoan2000} for more details on Kronecker product.
For the purpose of this paper, the most important property of Kronecker matrices is the possibility to efficiently compute matrix-vector products. 

Consider $d$ matrices $A^{(1)}, \ldots, A^{(d)}$, 
with $A^{(l)} \in \mathbb{R}^{s \times t}$, $l = 1,\ldots,d$. 
Let $(i_1,\ldots,i_d)$ and $(j_1,\ldots,j_d)$ be multi-indices with $1 \leq i_l \leq s $, $1 \leq j_l \leq t$, $l = 1,\ldots,d$. If we associate these multi-indices with the scalar indices 
\begin{equation}\label{eq:ordering} 
i  =  1 + \sum_{l=1}^d (i_l - 1) s^{l-1} \qquad \mbox{and} \qquad j = 1 + \sum_{l=1}^d (j_l - 1)t^{l-1},
\end{equation}
respectively, then from the definition of Kronecker product it follows
$$ \left(A^{(d)} \otimes \ldots  \otimes A^{(1)} \right)_{ij} = A^{(1)}_{i_1 j_1} \ldots A^{(d)}_{i_d j_d}. $$
Now assume we are given a vector $x \in \mathbb{R}^{t^d}$ and we want to compute the matrix-vector product $\left(A^{(d)} \otimes \cdots \otimes A^{(1)}\right) x$. We have
\begin{eqnarray*} \left(\left(A^{(d)} \otimes \ldots \otimes A^{(1)} \right) x \right)_i & = & \sum_{j = 1}^{t^d} \left(A^{(d)} \otimes \ldots \otimes A^{(1)} \right)_{ij} x_j \\
& = & \sum_{j_1,\ldots j_d = 1}^{t} A^{(1)}_{i_1 j_1} \ldots A^{(d)}_{i_d j_d} \mathcal{X}_{j_1, \ldots, j_d}, 
\end{eqnarray*}
where $\mathcal{X}_{j_1, \ldots, j_d}$ denotes the $(j_1, \ldots, j_d)$ entry of the $d-$dimensional tensor $\mathcal{X} \in \mathbb{R}^{t \times \ldots \times t}$, which is obtained by reshaping $x$ in accord with the second relation of \eqref{eq:ordering}.
If we reorder the terms in the last equality, we get the following sequence of nested summations
\begin{equation} \label{eq:sum-factorization}
\left(\left(A^{(d)} \otimes \ldots \otimes A^{(1)} \right) x \right)_i = \sum_{j_1 = 1}^{t} A^{(1)}_{i_1 j_1} \left( \sum_{j_2 = 1}^{t}  A^{(2)}_{i_2 j_2}\left( \ldots \sum_{j_d = 1}^{t} A^{(d)}_{i_d j_d} \mathcal{X}_{j_1, \ldots, j_d} \right)\right).
\end{equation}
By sequentially computing these nested sums, starting from the innermost one, for all values of the involved indices, the matrix-vector product can be computed at a low cost.
In the finite element/spectral element literature this procedure is often called \textit{sum-factorization}. In particular, by using \eqref{eq:sum-factorization}, the matrix $A^{(d)} \otimes \ldots \otimes A^{(1)}$ never needs to be formed explicitly. A complexity analysis reveals that the total number of flops needed to compute $\left(A^{(d)} \otimes \ldots \otimes A^{(1)} \right) x$ is  
$$ 2 \, st \left( t^{d-1} + s t^{d-2} + \ldots + s^{d-1} \right)$$
if the matrices $A^{(l)}$ are dense, and {\MT it is bounded by
\begin{equation} \label{eq:sparse_sumfact_flops} 2\, \max_{l} \mbox{nnz}\left(A^{(l)}\right)\left( t^{d-1} + s t^{d-2} + \ldots + s^{d-1} \right) \end{equation}
if the matrices $A^{(l)}$ are sparse, where $\mbox{nnz}(A^{(l)})$ denotes the number of nonzero entries of $A^{(l)}$.
}

\section{The MF-WQ approach} 
\label{sec:matrix-free}


When a linear system like \eqref{eq:linearsystem} is tackled with an iterative solver, e.g., the Conjugate Gradient (CG) method, the matrix $\A$ is required only to compute matrix-vector products. As a consequence, the formation and storage of the system matrix can be avoided entirely by providing a routine that takes in input a vector $v$ and returns the vector $\A v$ as output. This matrix-free approach can be very beneficial, both in terms of computational effort and memory storage. 

In this section we discuss how to use weighted quadrature in a  matrix-free approach for problem \eqref{eq:linearsystem}. 
We discuss how to compute matrix-vector products for  the mass matrix and the stiffness matrix, separately.

\subsection{The mass matrix}
\label{sec:mass}

Let $\widetilde{\M}$ be the approximation of $\M$ obtained with weighted quadrature, as described in Section \ref{sec:weighted-quadrature}. We want to compute the vector $\widetilde{\M} v$, where $v \in \mathbb{R}^N$ is a given vector.

For $i = 1,\ldots,N$, we observe that
$$ \left(\widetilde{\M} v\right)_i  = \sum_{j = 1}^N \widetilde{\M}_{ij} v_j = \sum_{j = 1}^N \sum_{q = 1}^{N_q} w_{i,q} c(\x_{q}) \hat{b}_j(\x_{q}) v_j = \sum_{q = 1}^{N_q} w_{i,q} c(\x_{q}) \left(\sum_{j = 1}^N \hat{b}_j(\x_{q}) v_j \right),$$
where in the second equality we used the definition of $\widetilde{\M}_{ij}$ from equation \eqref{eq:Mwq}. If we define
$$ v_h(\XI) = \sum_{j = 1}^N v_j \hat{b}_j (\XI), $$
we have then the obvious relation
\begin{equation} \label{eq:main}
\left(\widetilde{\M} v\right)_i = \sum_{q = 1}^{N_q} w_{i,q} \, c(\x_{q}) v_h(\x_{q}) = \Q_i \left( c(\cdot) v_h(\cdot) \right).  
\end{equation}
In other words, computing the $i-$th entry of $\widetilde{\M} v$ is equivalent to approximating the integral of the function $c \, v_h$ using the $i-$th quadrature rule.


MF-WQ  is based on equation \eqref{eq:main}. 
Indeed, it shows that 
$\widetilde{\M} v$
can be computed with the following steps:

\begin{enumerate}
\item Compute $\tilde{v} \in \mathbb{R}^{N_q}$, with $\tilde{v}_{q} := v_h(\x_{q})$, $q = 1,\ldots,N_q$.
\item Compute $\tilde{\tilde{v}} \in \mathbb{R}^{N_q}$, with $ \tilde{\tilde{v}}_{q} := c(\x_{q}) \cdot {\MT \tilde{v}_{q}}$, $q = 1,\ldots,N_q$.
\item Compute $\left(\widetilde{\M} v\right)_i = \sum_{q=1}^{N_q} w_{i,q} \, \tilde{\tilde{v}}_{q} $, $i = 1,\ldots,N$.
\end{enumerate}


This algorithm, and in particular steps 1 and 3, can be performed efficiently by exploiting the tensor structure of the basis functions and of the weights. In order to make this fact apparent, we now derive a matrix expression for the above algorithm. Consider the matrix of B-spline values $\mathcal{B} \in \mathbb{R}^{N_q \times N}$, with $\mathcal{B}_{q j} := \hat{b}_{j}(\x_{q})$, $q=1,\ldots,N_q$, $j = 1,\ldots,N$, which can be written as
\begin{equation} \label{eq:Bmass}  \mathcal{B} = \BB_d \otimes \ldots \otimes \BB_1, \end{equation}
where
\begin{equation} \label{eq:Bl} \left(\BB_l\right)_{q_l j_l} = \hat{b}_{l,j_l}(x_{l,q_l}), \qquad q_l = 1,\ldots,n_q, \; j_l = 1,\ldots,n. \end{equation}

We also consider the matrix of weights $\mathcal{W} \in \mathbb{R}^{N \times N_q}$, with $\mathcal{W}_{iq} := w_{i,q}$, $i = 1,\ldots, N$, $q=1,\ldots,N_q$. Thanks to the tensor structure of the weights, it holds
\begin{equation} \label{eq:Wmass} \mathcal{W} = \W_d \otimes \ldots \otimes \W_1, \end{equation}
where
$$ \left(\W_l\right)_{i_l q_l} = w_{l,i_l,q_l}, \qquad i_l = 1,\ldots,n, \; q_l = 1,\ldots,n_q. $$
Finally we introduce the diagonal matrix of coefficient values
\begin{equation} \label{eq:Dmass} \mathcal{D} := \mbox{diag}\left( \left\{c(\x_{q})\right\}_{q = 1,\ldots, N_q} \right). \end{equation}
Then for every $i,j = 1,\ldots,N$ from \eqref{eq:Mwq} we infer that
$$ \widetilde{\M}_{ij} = \sum_{q = 1}^{N_q} w_{i,q} c(\x_{q}) \hat{b}_{j}(\x_{q}) =  \sum_{q = 1}^{N_q} \mathcal{W}_{iq} \mathcal{D}_{qq} \mathcal{B}_{qj} = \left( \mathcal{W} \mathcal{D} \mathcal{B} \right)_{ij}. $$
Thus it holds
\begin{equation} \label{eq:WDB} \widetilde{\M} = \mathcal{W} \mathcal{D} \mathcal{B}  .\end{equation}
The factorization above of $\widetilde{\M}$ justifies Algorithm \ref{MF},  which computes efficiently the matrix-vector product.



\begin{algorithm} \LinesNumbered 
\SetKwInOut{Input}{Input}
\SetKwInOut{Output}{Output}
\SetKwInOut{Setup}{Setup} \Setup{Compute and store the matrices $\mathcal{D}$, $\BB_l$ and $\W_l$, for $l=1,\ldots,d$.}
\Input{ Vector $v \in \mathbb{R}^N$. }
Compute $\tilde{v} = \left(\BB_d \otimes \ldots \otimes \BB_1\right) v$  \;
Compute $\tilde{\tilde{v}} = \mathcal{D} \tilde{v}$ \;
Compute $w = \left(\W_d \otimes \ldots \otimes \W_1\right) \tilde{\tilde{v}} $ \;
\Output{ Vector $w = \widetilde{\M}v  \in \mathbb{R}^N$. }
\caption{MF-WQ product (mass)}\label{MF}
\end{algorithm}

We now analyze Algorithm \ref{MF} in terms of memory usage and of computational cost, where we distinguish between setup cost and application cost.
The Setup of Algorithm \ref{MF} requires the computation and storage of the coefficient values $c(\x_{q})$, $q=1,\ldots,N_q$, 
and of the (sparse) matrices $\W_l \in \mathbb{R}^{n\times n_q}$ and $\BB_l \in \mathbb{R}^{n_q\times n}$, for $l=1,\ldots,d$ {\MT (note that $\mathcal{W}$ and $\mathcal{B}$ never need to be formed)}.
The latter part, which involves only the computation and storage of univariate function values and weights, has negligible requirements both in terms of memory and arithmetic operations. 
The computational cost  of the evaluation of the coefficients $c(\x_{q})$ is problem dependent. 
{\MT
For our model problem, 
since the spline/NURBS parametrization $\F$ is fixed, 
and typically has a low degree (quadratics or cubics), the coefficients can be evaluated before the $k-$refinement is performed. 
Therefore, we can assume its cost is $O(N)$ flops,  independent of $p$.}
\\ {(secondo te dovrei dire che assumiamo che valutare $\alpha$ sia cheap??)}
%
In general, the storage of such coefficients clearly requires $N_q \approx 2^d N = O(N)$ memory
 (however, it is not necessary to store the whole $\mathcal{D} $, $ {\tilde{v}} $ and  $\tilde{\tilde{v}} $  since $w$ in Algorithm \ref{MF} can be computed component by component with on-the-fly calculation of the portion of $\mathcal{D} $, $ {\tilde{v}} $ and  $\tilde{\tilde{v}} $ that is needed,   see Remark \ref{rem:memory}).
We emphasize that this memory requirement is completely independent of $p$; this is a great improvement if we consider that storing the whole mass matrix would require roughly $(2p+1)^d N = O(N p^d)$ memory.
As for the application cost, Step 2 only requires $N$ flops. Using \eqref{eq:sparse_sumfact_flops} and the fact that $\mbox{nnz}\left(\BB^{(l)}\right) \approx 2pn$, $l = 1,\ldots,d$, we find that the number of flops required by Step 1 is
$$ 4pn \left(n^{d-1} + 2 n^{d-1} + \ldots + 2^{d-1} n^{d-1} \right) \leq 2^{d+2} N p = O(N p).  $$
Approximately the same number of operations is required for Step 3. Hence we conclude that the total application cost of Algorithm is $O(N p)$ flops. This should be compared with the $O(N p^d)$ cost of the standard matrix-vector product.

\subsection{The stiffness matrix} 
\label{sec:stiffness}

Matrix-vector products involving the stiffness matrix $\widetilde{\K}$ can be computed similarly as in the case of the mass matrix. Indeed, given $v \in \mathbb{R}^N$, by exploiting the definition of $\widetilde{\K}_{ij}$ from \eqref{eq:Kwq} we have that
$$ \left(\widetilde{\K} v\right)_i  =
\sum_{\alpha, \beta = 1}^d \sum_{j = 1}^N \sum_{q = 1}^{N_q} w^{(\alpha,\beta)}_{i,q} C_{\alpha \beta}(\x_{q}) \, \partial_{\beta} \hat{b}_j(\x_{q}) v_j = \sum_{\alpha, \beta = 1}^d \sum_{q = 1}^{N_q} w^{(\alpha,\beta)}_{i,q} C_{\alpha \beta}(\x_{q}) \,\partial_{\beta} v_h (\x_{q}),$$
where $v_h(\XI) = \sum_{j = 1}^{n} v_j \hat{b}_j(\XI)$.
Thus, $\widetilde{\K} v $ can be computed in the following steps:
\begin{enumerate}
\item Compute $\tilde{v}^{(\beta)} \in \mathbb{R}^{N_q}$, with $\left(\tilde{v}^{(\beta)}\right)_q := \partial_{\beta} v_h(\x_{q})$, $q = 1,\ldots,N_q$, for every $\beta = 1,\ldots,d$.
\item Compute $\tilde{\tilde{v}}^{(\alpha,\beta)} \in \mathbb{R}^{N_q}$, with $\left(\tilde{\tilde{v}}^{(\alpha,\beta)}\right)_q := C_{\alpha \beta}(\x_{q}) \cdot {\MT \left(\tilde{v}^{(\beta)}\right)_q}$, $q = 1,\ldots,N_q$, for every $\alpha,\beta = 1,\ldots,d$.
\item Compute $ w^{(\alpha,\beta)} \in \mathbb{R}^N$, with $ \left(w^{(\alpha,\beta)}\right)_i := \sum_{q = 1}^{N_q} w^{(\alpha,\beta)}_{i,q} \, \tilde{\tilde{v}}^{(\alpha,\beta)}_q $, $i = 1,\ldots,N$, for every $\alpha,\beta = 1,\ldots,d$.
\item Compute $\widetilde{\K} v = \sum_{\alpha,\beta = 1}^d w^{(\alpha,\beta)} $.
\end{enumerate}

Similarly as before, we derive a new matrix expression for $\widetilde{\K}$ which is the basis for MF-WQ.
For each $\beta = 1,\ldots,d$ we define $\mathcal{B}^{(\beta)} \in \mathbb{R}^{N_q \times N}$ as $\left( \mathcal{B}^{(\beta)}\right)_{qj} := \partial_{\beta} \hat{b}_{j}(\x_q) $, $q = 1,\ldots,N_q$, $j = 1, \ldots,N$. Equivalently, 
$$ \mathcal{B}^{(\beta)} = \BB^{\left(\beta\right)}_d \otimes \ldots \otimes \BB^{\left(\beta\right)}_1, $$
where
$$ \BB^{\left(\beta\right)}_l = \left\{\begin{array}{ll} \dot{\BB}_l & \mbox{if } l = q \\ \BB_l & \mbox{otherwise} \end{array}\right. $$
with $\BB_l$ defined as in \eqref{eq:Bl} and
$$ \left(\dot{\BB}_l\right)_{i_l q_l} = \hat{b}'_{l,j_l}(x_{l,q_l}), \qquad q_l = 1,\ldots,n_q, \; j_l = 1,\ldots,n. $$
Moreover, for every $\alpha, \beta = 1,\ldots,d$ we define $\mathcal{W}^{(\alpha,\beta)} \in \mathbb{R}^{N \times N_q}$ as $\left( \mathcal{W}^{(\alpha,\beta)}\right)_{iq} := w^{(\alpha,\beta)}_{i,q}$, $i = 1,\ldots,N$, $q = 1, \ldots,N_q$. Equivalently
$$ \mathcal{W}^{(\alpha,\beta)} = \W^{\left(\alpha,\beta\right)}_d \otimes \ldots \otimes \W^{\left(\alpha,\beta\right)}_1, $$
where
$$\left( \W^{\left(\alpha,\beta\right)}_d\right)_{i_l q_l} = w^{(\alpha,\beta)}_{l,i_l,q_l}, \qquad q_l = 1,\ldots,n_q, \; j_l = 1,\ldots,n.  $$
For $\alpha, \beta = 1,\ldots,d$ we finally define the diagonal matrix
$$ \left(\mathcal{D}^{\left(\alpha,\beta\right)}\right)_{qq} := C_{\alpha \beta}(\x_{q}), \qquad q = 1,\ldots,N_q. $$
Then from \eqref{eq:Kwq} by direct computation we find that
$$ \widetilde{\K} = \sum_{\alpha,\beta = 1}^d \mathcal{W}^{(\alpha,\beta)} \mathcal{D}^{\left(\alpha,\beta\right)} \mathcal{B}^{(\beta)}. $$
%
%
%
To compute matrix-vector products with $\widetilde{\K}$, we can hence use the following algorithm.

\begin{algorithm} \LinesNumbered 
\SetKwInOut{Input}{Input}
\SetKwInOut{Output}{Output}
\SetKwInOut{Setup}{Setup} \Setup{Compute and store the matrices $\mathcal{D}^{(\alpha,\beta)}$, $\W_l^{(\alpha,\beta)}$ and $\BB^{(\beta)}_l$ for $l,\alpha,\beta = 1,\ldots,d$.}
\Input{ Vector $v \in \mathbb{R}^N$. }
Initialize $w = 0$ \;
\For{$\beta=1,\dots, d$}
{Compute $\tilde{v}^{(\beta)} = \left(\BB_d^{(\beta)} \otimes \ldots \otimes \BB_1^{(\beta)} \right) v$ \;
\For{$\alpha=1,\dots, d$}{
Compute $\tilde{\tilde{v}}^{(\alpha,\beta)} = \mathcal{D}^{(\alpha,\beta)} \tilde{v}^{(\beta)}$ \;
Compute $ w^{(\alpha,\beta)} = \left( \W^{\left(\alpha,\beta\right)}_d \otimes \ldots \otimes \W^{\left(\alpha,\beta\right)}_1 \right) \tilde{\tilde{v}}^{(\alpha,\beta)} $ \;
Update $ w = w + w^{(\alpha,\beta)} $ \;
}
}
\Output{ Vector $w = \widetilde{\K}v  \in \mathbb{R}^N$. }
\caption{MF-WQ product (stiffness)}\label{MFstiff}
\end{algorithm}

As in the case of the mass, we analyze the memory and cost requirements of the MF-WQ product Algorithm \ref{MFstiff}.
The storage of the coefficients $C_{\alpha \beta}(\x_{q})$, which represents the main memory usage of Algorithm \ref{MFstiff}, requires $d^2$ times more memory than the mass algorithm, and precisely $d^2 N_q \approx d^2 2^d N = O(N) $ memory.
Note that this can be reduced to 
$\frac{d(d +1)}{2} N_q$ 
when, as in our example, symmetry  $C_{\alpha \beta} = C_{\beta \alpha}$ holds. See also Remark \ref{rem:memory} for a strategy that further reduces the memory consumption. 

Similarly as in the mass case, we can assume that the computation of such coefficients, which represents the main computational effort of the Setup step, require $O(N)$ flops.

As for the application cost, the presence of the two loops makes Algorithm \ref{MFstiff} roughly 
$ \frac{d(d +1)}{2} $ 
more expensive the its mass counterpart. We emphasize, however, that the cost is still $O(N p)$ flops.

\begin{remark} \label{rem:memory}
The memory required by Algorithm \ref{MFstiff} is $d^2 N_q$ (or $\frac{d \left(d+1\right)}{2}N_q$ by exploiting symmetry), due to the storage of the coefficient matrices $\mathcal{D}^{(\alpha,\beta)}$, $\alpha,\beta = 1,\ldots,d$.
However, we can consider a variant of Algorithm \ref{MFstiff} where these matrices are not stored and the coefficients $C_{\alpha \beta}(x_{q})$, $q=1,\ldots,N_q$, are computed on-the-fly every time $\tilde{\tilde{v}}^{(\alpha,\beta)}$ needs to be computed and then discarded.
In this variant, the main memory consumption is the storage of either vector $\tilde{\tilde{v}}^{(\alpha,\beta)}$ (note that, according to the structure of Algorithm \ref{MFstiff}, only one of these vectors, for $\alpha,\beta = 1,\ldots,d$, needs to be stored at a time), which requires $N_q$ memory.
Further memory saving can be achieved by computing  $w$  component by component and, for each component, on-the-fly calculation of the portion of $\mathcal{D}^{(\alpha,\beta)} $, $ {{\tilde{v}}^{(\beta)} } $ and  $\tilde{\tilde{v}}^{(\alpha,\beta)}  $ that is needed. This approach is not considered here. 
Of course, computing these coefficients on-the-fly increases the cost to compute matrix-vector products in terms of flops counting. 
{\MT Finally, we remark the same ideas can be applied to the mass case. However, the storage of $\mathcal{D}$ in Algorithm \ref{MFstiff} only requires $N_q$ memory, so it is probably less important to reduce memory consumption in this case.}
\end{remark}

\begin{remark}
In the derivation of MF-WQ,  the matrix entries stem from weighted quadrature.  However, Algorithms \ref{MF} and \ref{MFstiff} can be slightly modified to work with standard quadrature rules having a tensor structure. Consider for example the approximation of the mass matrix entries obtained with standard Gaussian quadrature
$$ \M^{GQ}_{ij} = \sum_{q } w^{GQ}_q c(\x_{q}^{GQ}) \hat{b}_{i}\left(\x_{q}^{GQ}\right) \hat{b}_{j}\left(\x_{q}^{GQ}\right) \qquad i,j = 1,\ldots,N, $$
where $\x_{q}^{GQ}$ and $w^{GQ}_q$ denote the Gaussian quadrature points and weights, respectively.
It can be easily checked that $ \M^{GQ}$ allows a matrix representation analogous to \eqref{eq:WDB}, i.e.
\begin{equation} \label{WDBgaussian} \M^{GQ} = \left(\W^{GQ}_d \otimes \ldots \otimes \W^{GQ}_1 \right) \mathcal{D}^{GQ} \left(\BB^{GQ}_d \otimes \ldots \otimes \BB^{GQ}_1 \right), \end{equation}
where $\mathcal{D}^{GQ}$ and $\BB^{GQ}_l$, $l = 1,\ldots,d$, are defined in analogy with \eqref{eq:Dmass} and \eqref{eq:Bl}, while $\W^{GQ}_l$, $l = 1,\ldots,d$, is the matrix of basis function values scaled by the weights, i.e.
$$ \left(\W^{GQ}_l\right)_{i_l q_l} = w^{GQ}_{l,q_l} \hat{b}_{l,i_l} (x_{l,q_l}). $$
However, the use of weighted quadrature presents a clear advantage, namely the reduced number of quadrature points. 
Indeed, Gaussian quadrature (and in general any quadrature rule that uses $O(p)$ points per knot span) requires $O(np)$ points for each direction. As a consequence, if \eqref{WDBgaussian} was used in a matrix-free algorithm analogous to Algorithm \ref{MF}, the computation and storage of $\mathcal{D}^{GQ}$ in the Setup step would require $O(N p^d)$ flops and memory, rather than just $O(N)$. Moreover, the main computational effort of each matrix-vector product, namely the computation of $\left(\BB^{GQ}_d \otimes \ldots \otimes \BB^{GQ}_1\right) v$ and $\left(\W^{GQ}_d \otimes \ldots \otimes \W^{GQ}_1\right) \tilde{\tilde{v}}$, would require $O(N p^{d+1})$ flops, rather than just $O(N p)$. The same observation applies  to the stiffness matrix case.
\end{remark}

\section{Numerical experiments} 
\label{sec:experiments}

\begin{figure}
\begin{center}
\includegraphics[trim=0cm 0cm 0cm 0cm, clip=true, width=0.5\textwidth]{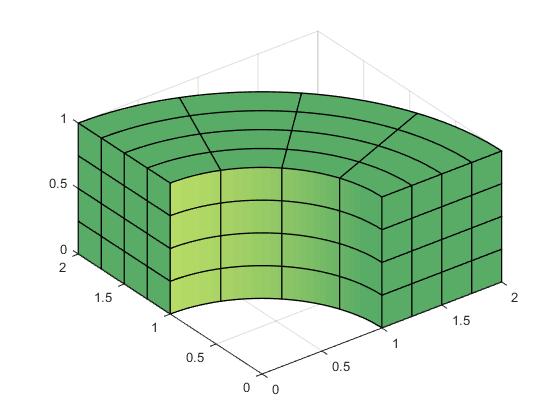}
\caption{Thick ring domain}
\label{fig:thickring}
\end{center}
\end{figure}

In this section, we show numerical evidence of the potential of the MF-WQ approach, analyzing both  memory and CPU usage. The focus is on three-dimensional problems. Accordingly, in all our experiments we take as domain the thick quarter of ring 
$$ \Omega = \left\{ (x_1,x_2,x_3) \in \mathbb{R}^3 \left.\right| \, 1 \leq x_1^2 + x_2^2 \leq 4, \, 0 \leq x_1,x_2 , \, 0 \leq x_3 \leq 1 \right\},$$
which is represented in Figure \ref{fig:thickring}.


All the experiments presented here are performed in Matlab, version 8.5.0.197613 (R2015a), on a  
Linux workstation equipped with Intel i7-5820K processors running at 3.30GHz, and with 64 GB of RAM.
Since parallelization is not considered in this work we benchmark only sequential execution times, and consistently all the computations are restricted to one computational thread.

\begin{figure}
\begin{center}
 \includegraphics[width=0.6\textwidth] {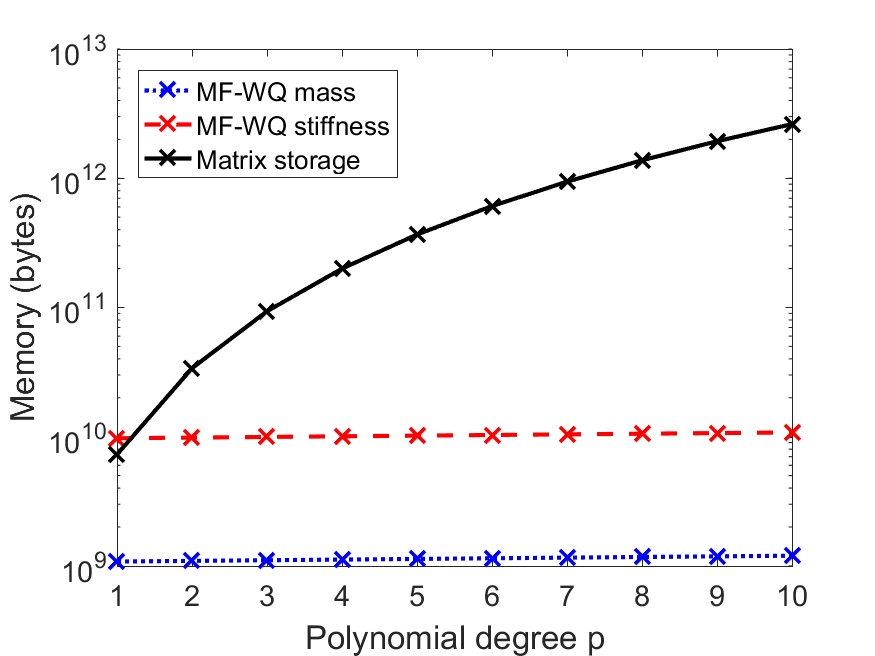} 
\end{center}
\caption{ Memory consumption for the MF-WQ approach for the mass and the stiffness, and for the matrix storage. } 
\label{fig:memory}
\end{figure}

We start by considering the memory requirements, which in MF-WQ corresponds to the memory required to store the matrices $\mathcal{D}$, $\W_l$, $\BB_l$, $l = 1,2,3$, in the mass case, and $\mathcal{D}^{\left(\alpha,\beta\right)}$, $\W_l^{\left(\alpha,\beta\right)}$, $\BB_l^{\left(\beta\right)}$, $l,\alpha,\beta = 1,2,3$, in the stiffness case. In Figure \ref{fig:memory}, we display the memory occupied in both cases, for a uniform discretization of $\Omega$ with { $256^3$} elements, and for different values of $p$. We remark that, for ease of implementation, we did not take advantage of the fact that $C_{\alpha \beta} = C_{\beta \alpha}$ (and hence $\mathcal{D}^{\left(\alpha,\beta\right)} = \mathcal{D}^{\left(\beta,\alpha\right)}$ ) and did not consider the memory saving techniques described in Remark \ref{rem:memory}. Indeed, these are not necessary in our tests. As a comparison, in Figure \ref{fig:memory} we also display the memory that would be required to actually store either matrix\footnote{Since for large enough $p$ it would be impossible to store matrices with our RAM resources, the memory reported in Figure \ref{fig:memory} is actually just an approximation. This is given by 16 bytes for each nonzero entries (the number of nonzero entries is computed exactly). This is indeed a lower bound, and also a good approximation, for the memory occupied by a sparse matrix in the 64-bit version of Matlab.} 
(note that both matrices require the same memory, since they have exactly the same dimension and the same sparsity pattern). 
We emphasize that this is the minimum possible memory required to form the matrix. Depending on the actual implementation, the formation algorithm might require an even greater amount of memory.

{ As it can be seen in the plot, MF-WQ  for the mass requires about 1 GB, while for the stiffness about 10 GB. We emphasize that, as expected, the memory consumption is almost constant with respect to $p$. This is of course not true when the matrix is stored. In this case, the memory consumption exceeds 64 GB already for $p=3$, 256 GB for $p = 5$, and 1 TB for $p=8$. }

\begin{figure}
\begin{center}
 \includegraphics[width=0.49\textwidth] {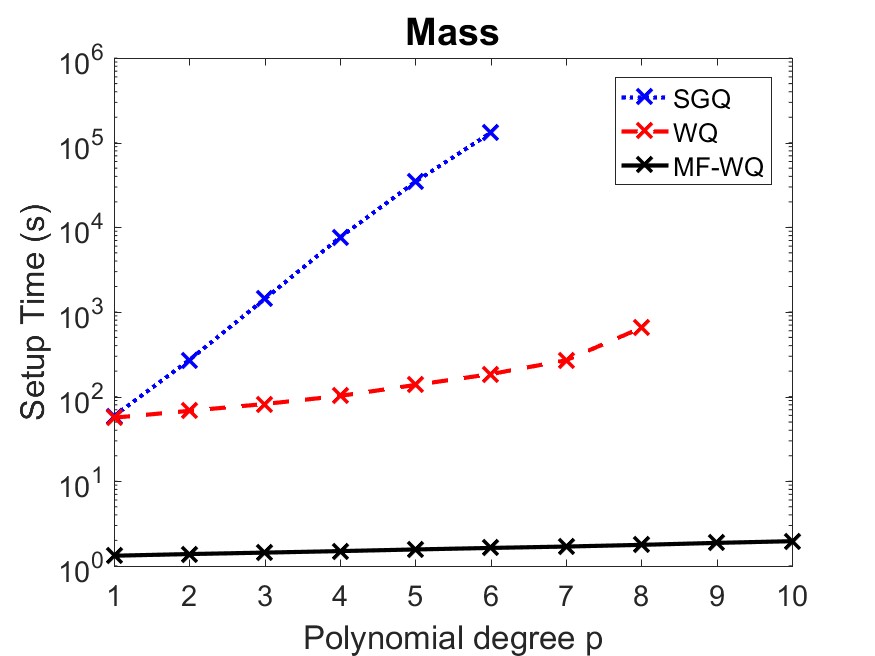} 
 \includegraphics[width=0.49\textwidth] {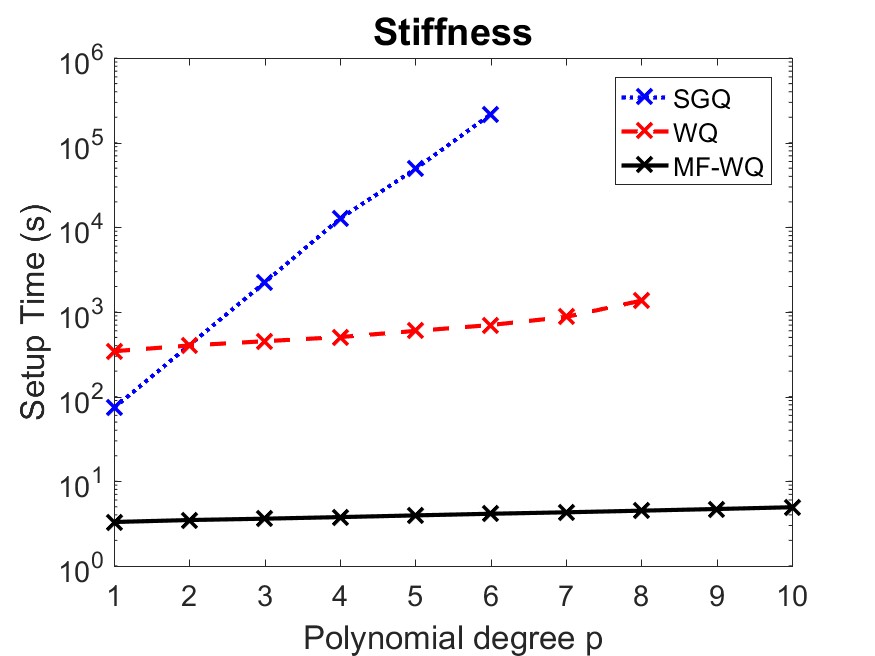} 
\end{center}
\caption{Setup time for the mass (left plot) and for the stiffness (right plot), when the matrix is formed with SGQ and with WQ, and when using the MF-WQ approach.} 
\label{fig:setup_time}
\end{figure}

We now benchmark the computation time required for the Setup of Algorithms \ref{MF} and \ref{MFstiff}, when using a uniform discretization on $\Omega$ with $64^3$ elements. In this case, the performance is compared with that of two approaches where the Galerkin matrix is constructed  and stored: in this case we consider both standard Gaussian quadrature (SGQ), as implemented by the Matlab toolbox GeoPDEs 3.0 \cite{vazquez2016new}, and weighted quadrature (WQ), as in  in \cite{Calabro2017}. 
Times are computed by the {\tt tic} and {\tt toc} commands of Matlab. 
Note that we are considering a coarser discretization than in the previous example, since in this way we can afford storing the matrices even for rather large values of $p$. Results are reported in Figure \ref{fig:setup_time}, where the left plot refers to the mass matrix, while the right plot refers to the stiffness matrix.

As it can be seen from the plot, the difference between  MF-WQ and the other approaches is very noticeable. Indeed, even for $p=1$ the matrix-free setup requires between 1 and 2 orders of magnitude less time than the fastest matrix formation algorithm, both for the mass and for the stiffness; and clearly the gap becomes wider as $p$ is increased. In particular, as expected the setup time required by  MF-WQ  is almost constant with respect to $p$.

\begin{figure}
\begin{center}
 \includegraphics[width=0.6\textwidth] {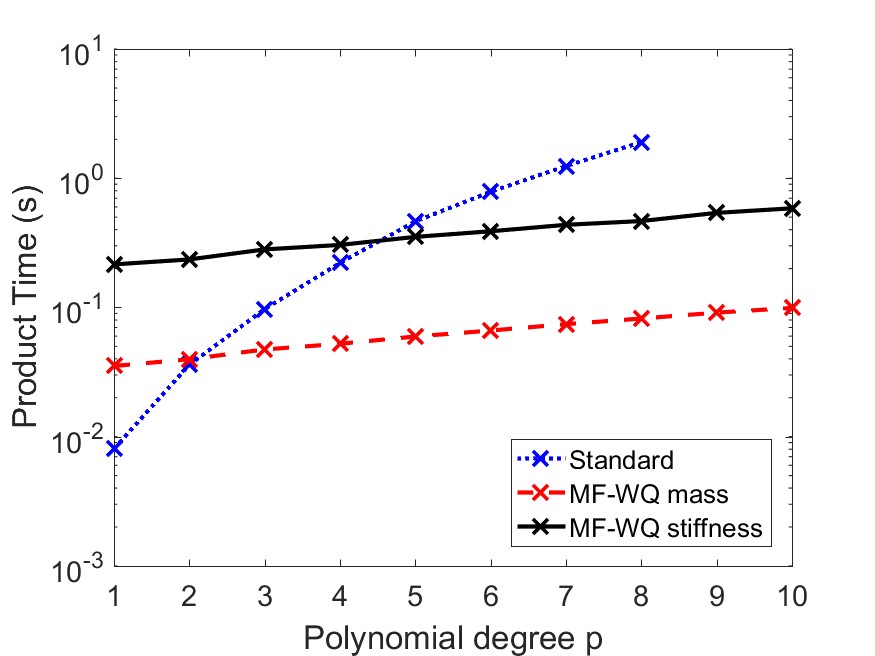} 
\end{center}
\caption{Computation time for a single matrix-vector product with MF-WQ for the mass and for the stiffness, and with standard product of Matlab.} 
\label{fig:product_time}
\end{figure}

For the same discretization, in Figure \ref{fig:product_time} we plot the time spent to compute a single matrix-vector product with Algorithms \ref{MF} and \ref{MFstiff}, and with the standard product of Matlab when the matrix is stored. 
{ \MT We remark that this comparison somewhat penalizes our approach, since the Matlab product is a built-in operation while our strategy relies on interpreted code.
Nevertheless, the results show that MF-WQ  product is faster, except for lower degrees.
Note also that MF-WQ  product includes quadrature, which is part  of  the matrix construction (setup) for the standard product. In any case, setup time likely dominates the matrix-vector product time (see the next set of experiments). }

Up to this point, we have only analyzed single aspects of MF-WQ.  Now we test its overall efficiency for solving a  differential model problem. 
We consider a Dirichlet problem of the form \eqref{eq:poisson}, {\MT with $\alpha = 0$ and $K = I$}, on a mesh of $256^3$ elements. For the sake of simplicity, we consider a uniform mesh  but all the algorithms considered do not take any advantage of it  and  are written to work on non-uniform meshes.

In order to motivate  the use of high-order splines on a fine mesh, we consider an oscillating manufactured solution, namely 
\begin{equation}
  \label{eq:exact-solution}
  \mathrm{u} (x_1,x_2,x_3) = \sin\left(5 \pi x_1\right) \sin\left(5 \pi x_2\right) \sin\left(5 \pi x_3\right) \left(x_1^2 + x_2^2 - 1\right)\left(x_1^2 + x_2^2 - 4\right). 
\end{equation}
An oscillating solution is somewhat unnatural for our model  problem, but it is representative of solutions that may arise in more complicated differential problems (such as  Navier-Stokes in the turbulent regime). 

In addition to MF-WQ, we consider forming the linear systems by SGQ and WQ.
In each case,  we have to select a method to solve the linear system. In MF-WQ, an iterative method is the only viable option. On the other hand, when the system matrix is stored either an iterative or a direct method can be used.
However, it is well-known that direct solvers usually require a huge amount of memory, which might well exceed the available memory.
Indeed, in our experiments Matlab's sparse direct solver ``backslash'' fails to solve the system even for $p = 1$. Thus, we rely on an iterative approach also when the system matrix is stored.

When we use SGQ,  we take advantage of the symmetry of the system matrix and  solve the system by the Conjugate Gradient (CG) method. On the other hand, for WQ and  MF-WQ, since the matrix is nonsymmetric, we use BiCGStab. 

Computational efficiency of iterative solvers requires a careful choice of the stopping criterion. The rule of thumb is to solve the linear system up to an accuracy that matches the Galerkin error 
{ $\left\|\mathrm{u}_h - \mathrm{u}\right\|_{H^1}$, where $\mathrm{u}_h$ is the {\MT (exact)} Galerkin solution.
This motivates our choice for the stopping tolerance on the residual as follows:
\begin{equation} \label{eq:stopping_criterion}\frac{\left\| \textbf{r}_k \right\|_2}{\left\|\f\right\|_2} \leq \eta \frac{\left\|\mathrm{u}_h - \mathrm{u}\right\|_{H^1}}{\left\|\mathrm{u}\right\|_{H^1}} ,\end{equation}
where $\textbf{r}_k$ is the residual at the $k-$th iteration, and with $\eta = 0.1$. In our case,}  the relative Galerkin error at the right hand side is estimated  knowing the exact solution \eqref{eq:exact-solution} (In real applications, the Galerkin error is unknown  and needs to be  approximated by a suitable a posteriori error estimator).  

{\MT
As a preconditioner, we consider the Fast Diagonalization (FD) method. This approach, originally proposed in \cite{Lynch1964}, was first used in the context of isogeometric analysis in \cite{Sangalli2016}. \GS The FD preconditioner has been further developed in \cite{montardini2017robust} in order to improve its performance on complex geometries. However, for the sake of simplicity here we use its simplest version from \cite{Sangalli2016}. \B In \cite{Tani2017} the FD preconditioner is discussed in the context of weighted quadrature, assessing its robustness with respect to both the mesh size and the spline degree. 
}

\begin{figure}
\begin{center}
\includegraphics[width=0.6\textwidth] {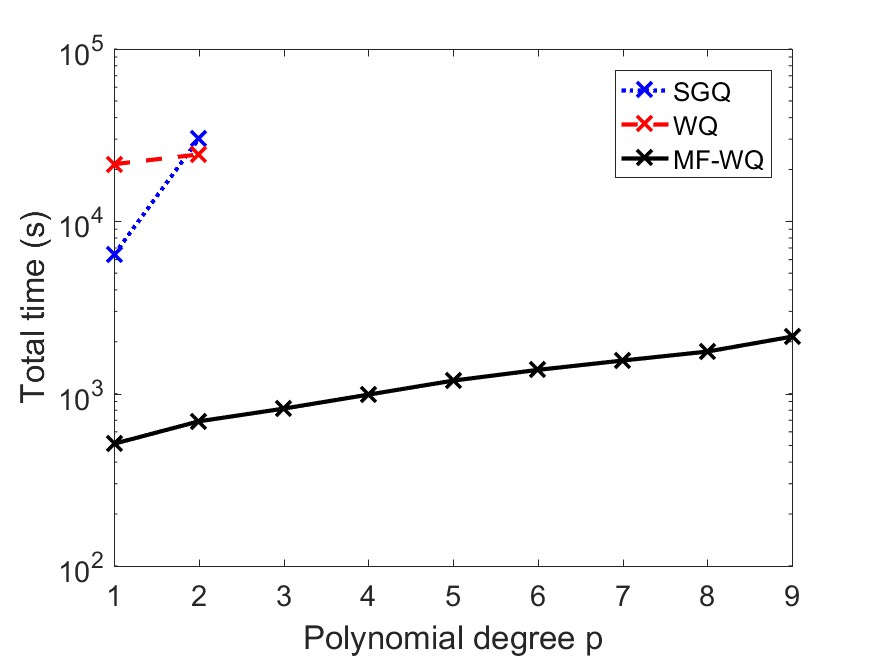} 
\end{center}
\caption{Total computation time for the solution of the Dirichlet problem when forming the matrix with the  the SGQ, WQ  and MF-WQ approaches.} 
\label{fig:diff_problem}
\end{figure}

In Figure \ref{fig:diff_problem} we report, for different values of $p$, the total computation time for the solution of the Dirichlet problem, which gathers the time spent for the setup of the linear system and the time spent for its solution.
Coherently with the results of Figure \ref{fig:memory}, the storage of the system matrix in the SGQ approach exceeds our memory limitations for $p \geq 3$, while in MF-WQ  all the tested values of $p$ are allowed.
Moreover, even for $p=1,2$ MF-WQ  appears to be much faster. Indeed, the speedup factor is $13$ times for $p=1$ and
44 times for $p=2$.

\begin{figure}
\begin{center}
\includegraphics[width=0.6\textwidth] {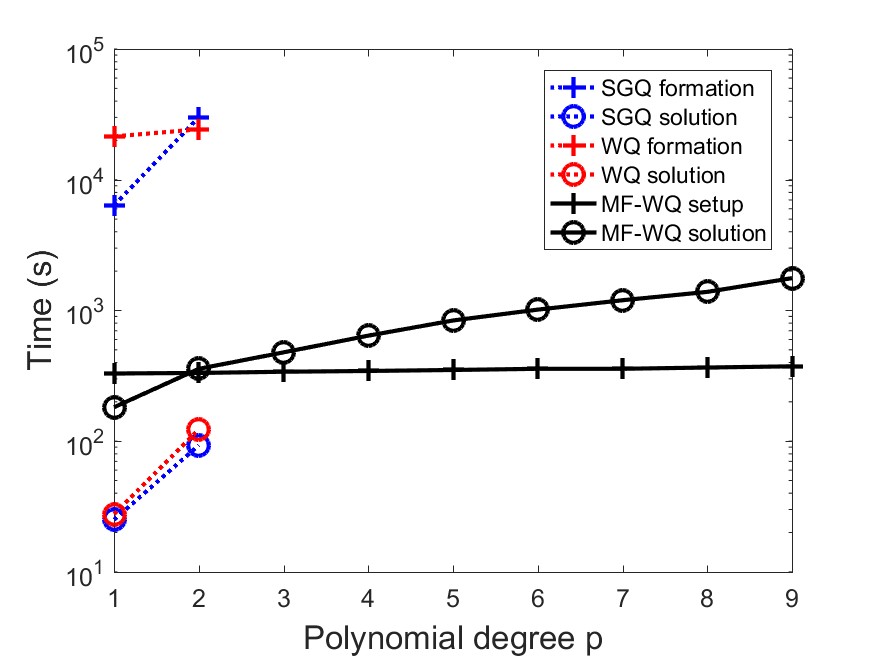} 
\end{center}
\caption{System setup time and system solution time for the Dirichlet problem when forming the matrix with the SGQ, WQ  and MF-WQ approaches.} 
\label{fig:diff_problem2}
\end{figure}

In Figure \ref{fig:diff_problem2} we report the same results, but splitting between setup time and solution time. 
It is interesting to observe that when the matrix is formed with SGQ the setup time is clearly dominant, while in MF-WQ setup and solution times are more balanced.


\begin{figure}
\begin{center}
\includegraphics[width=0.8\textwidth]{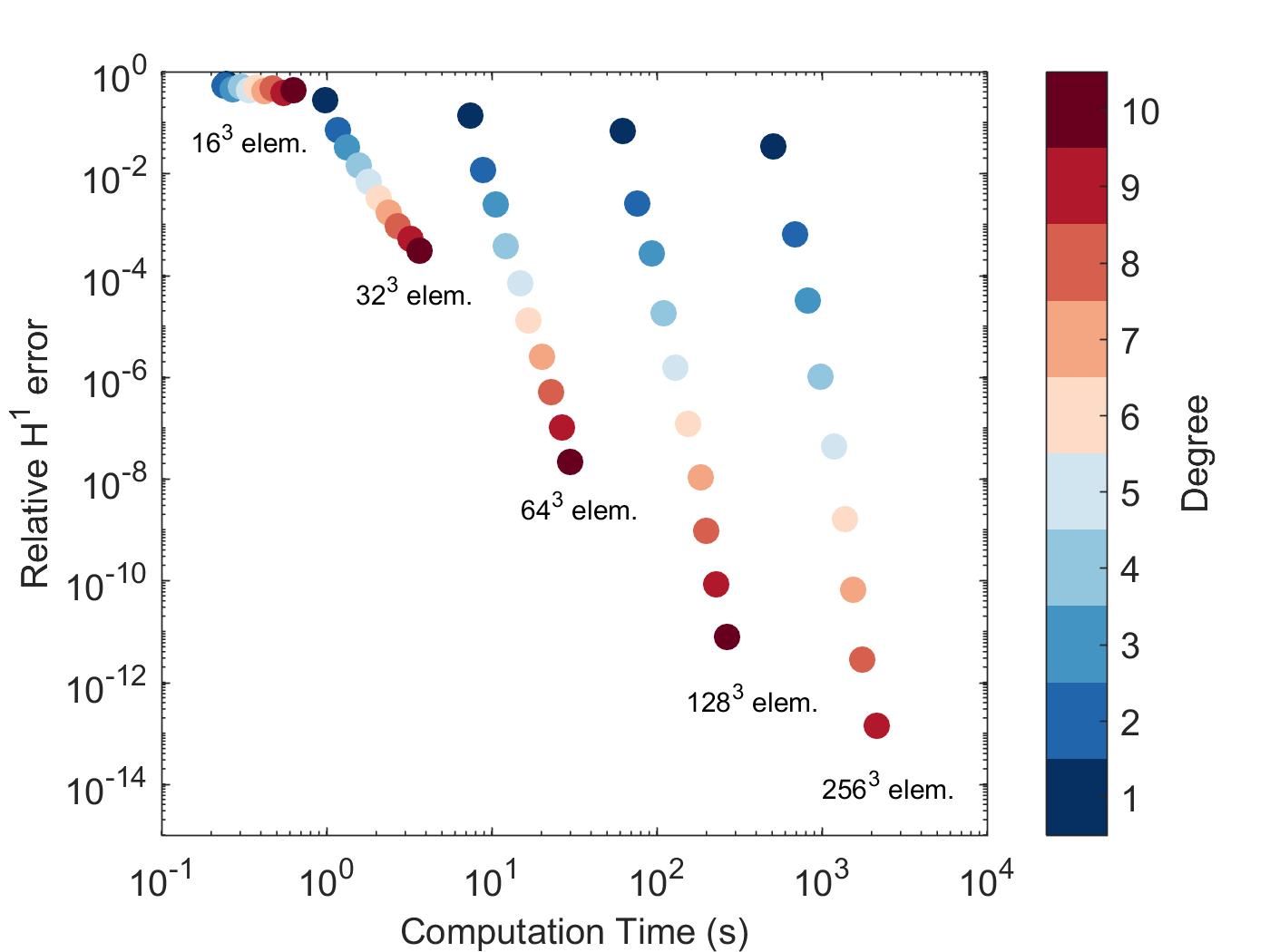} 
\end{center}
\caption{Representation in the time-error plane of MF-WQ results shown in Table \ref{tab:k_refinement}.} 
\label{fig:k_refinement}
\end{figure}

Our last set of experiments  aim at showing the most relevant result of this paper:  the isogeometric $k-$refinement, whose use has always been discouraged by its prohibitive computational cost,  becomes very appealing in the present setting. We still use as a benchmark the Poisson problem with solution \eqref{eq:exact-solution},  and solve it for different parametric mesh-size $h$ reaching  $h=2^{-8}$ (that is, the finest mesh is formed by $256^3$ elements)  and degree {\MT$p=1, \ldots, 9,10$}. We only consider MF-WQ, with the Krylov method, the preconditioner and the stopping criterion as above.
For each value of $h$ and $p$ we measure the total computation time (setup and solution of the system) and the error
$\left\|\mathrm{u} - \tilde{\mathrm{u}}_h \right\|_{H^1} / \left\|\mathrm{u} \right\|_{H^1}$, where $\tilde{\mathrm{u}}_h \in V_h$ is the function associated with the approximate solution of the linear system. Results are shown in Table \ref{tab:k_refinement} and in Figure \ref{fig:k_refinement} {\MT (in the latter, only for $h=2^{-8}$ the case $p=10$ is not reported, since the error reaches machine precision already for $p=9$, and precisely $\left\|\mathrm{u} - \tilde{\mathrm{u}}_h \right\|_{L^2} / \left\|\mathrm{u} \right\|_{L^2} \approx 10^{-15}$)}. We have the following observations.
\begin{itemize}
\item There is a minimal mesh resolution which is required to allow  $k$-refinement convergence.  This depends on the solution, which is in the example \eqref{eq:exact-solution}  a simple oscillating function  with wavelength $1/5$ on a domain with diameter $3$.  Indeed, there is no approximation (i.e., the relative approximating error remains close to $1$) for meshes of  $16^3$ elements or coarser, for any $p$. Convergence begin at a resolution of $32^3$ elements.
\item 

The computation time of MF-WQ  grows almost linearly with respect to $N = \left ( \frac{1}{h} \right )^{3}$ (note that the growth is slower between the two coarser discretization level, where apparently we are still in the pre-asymptotic regime).
Time dependence on $p$ is also very mild: {\MT the computation times for $p=10$ are only roughly 3 times larger than for $p=2$, keeping the same mesh resolution.}
We remark that the time growth with respect to $N$ and $p$ is due not only to the increased cost for system setup, matrix-vector product and application of the preconditioner, but also to the increased number of iterations. In turn, the number of iterations grows 
not because of a worsening of the preconditioner's quality (according to the results in \cite{Sangalli2016,Tani2017}) but because of a smaller discretization error, which corresponds to a more stringent stopping criterion according to \eqref{eq:stopping_criterion}.


\item The higher  the degree,  the higher  the computational efficiency of the $k$-method. This is clearly seen in Figure \ref{fig:k_refinement} where the red dots (associated to $p=9,10$, the highest degree in our experiments), are at the bottom of the error vs. computation time plot.  
\item The $k$-refinement is superior to low-degree $h$-refinement given a target accuracy:  for example, for a target accuracy of order $10^{-3}$, we can select degree $p=8$ on a mesh of $32^3$ elements  or $p=2$  on a mesh of $256^3$ elements: the former approximation is obtained in {$2.7$ seconds while the latter takes about $690$ seconds on our workstation, with  speedup factor higher than $250$.}
\end{itemize}

  \begin{remark}
    In the tests of this section we have not encountered any numerical instability, having considered degree up to $p=10$, that significantly exceeds the degrees used  in the isogeometric literature. 
	This is an indication that standard spline routines are robust and the spline basis,  whose condition number is known to depend exponentially from  $p$, in practice is well suited  for this range of degrees. It is possible that instabilities arise when higher polynomial degrees are considered, and further study would be needed to gain insight on this topic. 
  \end{remark}
	
{\MT
\begin{sidewaystable}
    \centering
\renewcommand\arraystretch{1.4} 
\begin{tabular}{|c|c|c|c|c|c|}
\hline
& \multicolumn{5}{|c|}{Relative $H^1$ error / Total computation time (s)} \\
\hline
\backslashbox{$p$}{$h$} & $2^{-4}$ & $2^{-5}$ & $2^{-6}$ & $2^{-7}$ & $2^{-8}$ \\
\hline
1 &  $5.8\cdot 10^{-1}$ / $2.4\cdot 10^{-1}$ & $2.8\cdot 10^{-1}$ / $9.9\cdot 10^{-1}$ & $1.4\cdot 10^{-1}$ / $7.4\cdot 10^{0}$ & $6.8\cdot 10^{-2}$ / $6.2\cdot 10^{1}$ & $3.4\cdot 10^{-2}$ / $5.1\cdot 10^{2}$ \\
\hline
2 & $5.3\cdot 10^{-1}$ / $2.5\cdot 10^{-1}$ & $7.1\cdot 10^{-2}$ / $1.1\cdot 10^{0}$ & $1.2\cdot 10^{-2}$ / $8.8\cdot 10^{0}$ & $2.6\cdot 10^{-3}$ / $7.6\cdot 10^{1}$ & $6.2\cdot 10^{-4}$ / $6.9\cdot 10^{2}$ \\
\hline
3 & $4.5\cdot 10^{-1}$ / $2.7\cdot 10^{-1}$ & $3.3\cdot 10^{-2}$ / $1.3\cdot 10^{0}$ & $2.5\cdot 10^{-3}$ / $1.1\cdot 10^{1}$ & $2.7\cdot 10^{-4}$ / $9.3\cdot 10^{1}$ & $3.2\cdot 10^{-5}$ / $8.2\cdot 10^{2}$ \\
\hline
4 & $5.1\cdot 10^{-1}$ / $3.0\cdot 10^{-1}$ & $1.4\cdot 10^{-2}$ / $1.6\cdot 10^{0}$ & $3.8\cdot 10^{-4}$ / $1.2\cdot 10^{1}$ & $1.8\cdot 10^{-5}$ / $1.1\cdot 10^{2}$ & $1.0\cdot 10^{-6}$ / $9.9\cdot 10^{2}$ \\
\hline
5 & $4.4\cdot 10^{-1}$ / $3.4\cdot 10^{-1}$ & $6.8\cdot 10^{-3}$ / $1.8\cdot 10^{0}$ & $7.1\cdot 10^{-5}$ / $1.5\cdot 10^{1}$ & $1.5\cdot 10^{-6}$ / $1.3\cdot 10^{2}$ & $4.3\cdot 10^{-8}$ / $1.2\cdot 10^{3}$ \\
\hline
6 & $4.9\cdot 10^{-1}$ / $3.8\cdot 10^{-1}$ & $3.3\cdot 10^{-2}$ / $2.1\cdot 10^{0}$ & $1.3\cdot 10^{-5}$ / $1.7\cdot 10^{1}$ & $1.2\cdot 10^{-7}$ / $1.6\cdot 10^{2}$ & $1.6\cdot 10^{-9}$ / $1.4\cdot 10^{3}$ \\
\hline
7 & $4.1\cdot 10^{-1}$ / $4.2\cdot 10^{-1}$ & $1.7\cdot 10^{-3}$ / $2.4\cdot 10^{0}$ & $2.5\cdot 10^{-6}$ / $2.0\cdot 10^{1}$ & $1.1\cdot 10^{-8}$ / $1.9\cdot 10^{2}$ & $6.7\cdot 10^{-11}$ / $1.6\cdot 10^{3}$ \\
\hline
8 & $4.7\cdot 10^{-1}$ / $4.7\cdot 10^{-1}$ & $9.2\cdot 10^{-4}$ / $2.7\cdot 10^{0}$ & $5.1\cdot 10^{-7}$ / $2.3\cdot 10^{1}$ & $9.3\cdot 10^{-10}$ / $2.0\cdot 10^{2}$ & $2.8\cdot 10^{-12}$ / $1.8\cdot 10^{3}$ \\
\hline
 9 & $3.8\cdot 10^{-1}$ / $5.5\cdot 10^{-1}$ & $5.2\cdot 10^{-4}$ / $3.2\cdot 10^{0}$ & $1.0\cdot 10^{-7}$ / $2.7\cdot 10^{1}$ & $8.4\cdot 10^{-11}$ / $2.3\cdot 10^{2}$ & $1.4\cdot 10^{-13}$ / $2.1\cdot 10^{3}$ \\
\hline
10 & $4.4\cdot 10^{-1}$ / $6.3\cdot 10^{-1}$ & $3.0\cdot 10^{-4}$ / $3.7\cdot 10^{0}$ & $2.2\cdot 10^{-8}$ / $3.0\cdot 10^{1}$ & $7.8\cdot 10^{-12}$ / $2.7\cdot 10^{2}$ & $2.8\cdot 10^{-13}$ / $2.0\cdot 10^{3}$ \\
\hline
\end{tabular} 
  \caption{$H^1$ error and total computation time for MF-WQ.}
\label{tab:k_refinement}
\end{sidewaystable}
}

\section{Conclusions} 
\label{sec:conclusions}
\GS 
We have proposed MF-WQ, a matrix-free approach for isogeometric analysis based on weighted quadrature.
Our goal has been  to reduce the huge computational effort associated with the formation of the Galerkin matrices, especially for the high-degree $k$-method, that is,  when splines of high-degree and high-continuity are adopted.

The memory required by MF-WQ is proportional to the number of  degrees-of-freedom of the discretization  and independent of $p$ when considering $k$-refinement.  Moreover the computational cost grows only linearly with respect to the degree $p$  and is orders of magnitude less than the cost of the standard formation of the matrix, based on element-wise Gaussian quadrature.

We test MF-WQ quadrature on a Poisson model problem, for which we have at our disposal an efficient preconditioner based on the Fast Diagonalization method, that we have proposed in a previous work  \cite{Sangalli2016}. Furthermore we target a case with a smooth solution. These are the ideal conditions for the  $k$-method,  which allow us to show the full  benefit of  MF-WQ: increasing the degree  and continuity leads to orders of magnitude higher computational efficiency.
This is shown for the first time, to our knowledge, in numerical experiments.

Further researches  will focus  on the extension of MF-WQ  to locally-refinable, trimmed and multi-patch  isogeometric spaces. Locally-refinable splines are needed  to  optimally approximate solutions with singularities, while trimmed and multi-patch geometries are used to represent complex  domains of interest. The multi-patch case is trivial to treat, the rest requires some re-thinking of the weighted quadrature idea.

The extension of  the  overall solver beyond the Poisson model problem considered in this paper also requires an efficient preconditioner. We remark that the FD preconditioner has been generalized in  \cite{montardini2017robust}, that  considers the Stokes system and improves the preconditioner robustness with respect to the geometry parametrization. Furthermore the use of the FD preconditioner for conforming  multi-patch domains is discussed in \cite{Sangalli2016}. However, there are important  challenges to face for the preconditioner development on realistic  geometries (e.g., trimmed  domains),  locally-refinable spaces, and more general  differential operators. The latter difficulty  is common to all numerical methods, in particular high-order ones. On the other hand, handling trimming or complex  geometries  (e.g., singular or highly distorted) is a specific difficulty of the isogeometric method,  and it is under study in the community. The same is true for  locally-refinable spaces, see the recent paper  \cite{cho2017bpx} and the references therein.
\B

\section*{Acknowledgments}

The authors were partially supported by the European Research Council
through the FP7 Ideas Consolidator Grant \emph{HIGEOM} n.616563.
The authors are members of the  Gruppo Nazionale Calcolo
Scientifico-Istituto Nazionale di Alta Matematica (GNCS-INDAM), and
the second  author was partially supported by GNCS-INDAM for this research.
This support is gratefully acknowledged.

\bibliography{biblio}

\end{document}